\documentclass[12pt,reqno]{amsart}
\usepackage{amssymb,amsmath,dsfont}
\usepackage{pgfplots}
\pgfplotsset{compat=1.14}
\usepackage{enumitem}
\usepackage{calc}
\usepackage{mathrsfs}
\usepackage[percent]{overpic}

\title[No coexistence in the two-type Richardson model]{No coexistence in the two-type Richardson model when one side is stronger than the other}
\date{}

\author{Gideon Amir}
\author{Dor Elboim}
\author{Gady Kozma}

\usepackage[margin=1in]{geometry}
\usepackage[all]{xy}
\usepackage{amssymb}
\usepackage{amsfonts}
\usepackage{amsthm}
\usepackage{amsmath}
\usepackage{mathtools}
\usepackage{url}
\usepackage{xcolor}
\usepackage{dsfont}
\usepackage{pgfplots}
\usepackage{bbm}
\usepackage{hyperref}

\definecolor{myblue}{rgb}{0.09,0.32,0.44} 
\hypersetup{pdfborder={0 0 0},pdfborderstyle={},colorlinks=true,linkcolor=myblue,citecolor=myblue,urlcolor=blue}

\newtheorem{thm}{Theorem}
\newtheorem*{thm*}{Theorem}

\newtheorem{lem}{Lemma}[section]  
\newtheorem{prop}[lem]{Proposition}

\newtheorem{claim}[lem]{Claim}

\theoremstyle{remark}

\newtheorem{remark}[lem]{Remark}

\newcommand{\ZZ}{\mathbb{Z}}

\makeatletter
\ifdefined\andify
\renewcommand{\andify}{%
  \nxandlist{\unskip, }{\unskip{} \@@and~}{\unskip{} \@@and~}}
\def\author@andify{%
  \nxandlist {\unskip ,\penalty-1 \space\ignorespaces}%
    {\unskip {} \@@and~}%
    {\unskip \penalty-2 \space \@@and~}%
}
\fi
\let\@wraptoccontribs\wraptoccontribs
\makeatother

\numberwithin{equation}{section}

\begin{document}

\begin{abstract}
    We study the two-type Richardson model in two dimensions. We prove that there is no coexistence of the two species when their infection rates are different, establishing a conjecture of H{\"a}ggstr{\"o}m and Pemantle from 1997.

    In our proof we analyze a second order approximation of the infected region and use a Mermin--Wagner style argument to push the strong type out of the weak one. 
\end{abstract}

\maketitle

\begin{figure}[htbp]
    \centering
    \resizebox{\linewidth}{!}{%
        \includegraphics[height=3cm]{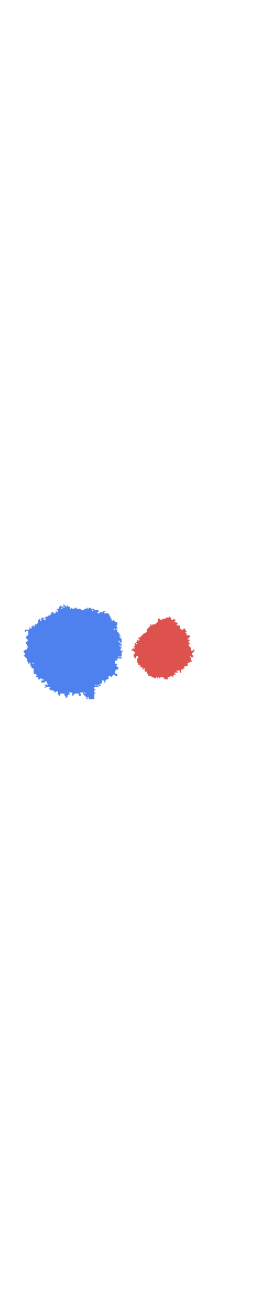}\hspace{2mm}%
        \includegraphics[height=3cm]{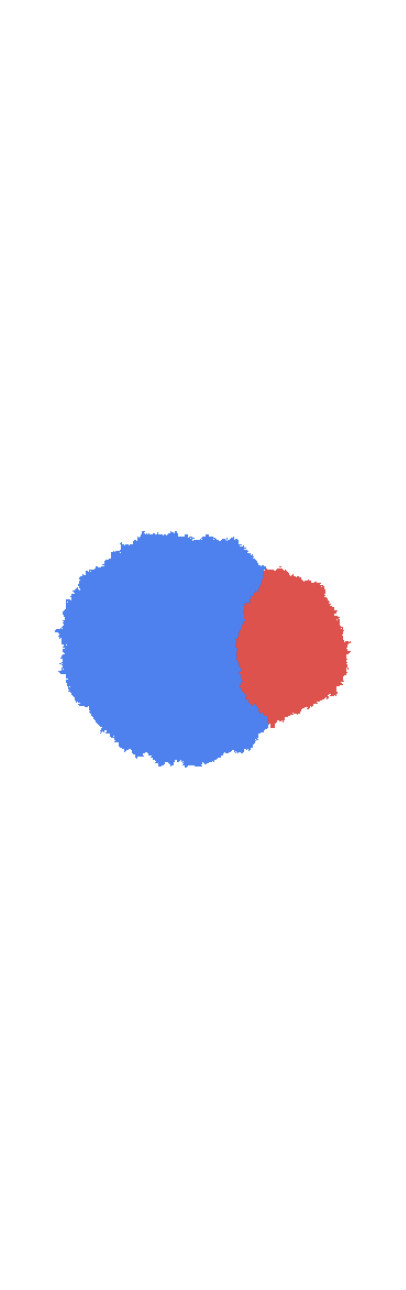}\hspace{2mm}%
        \includegraphics[height=3cm]{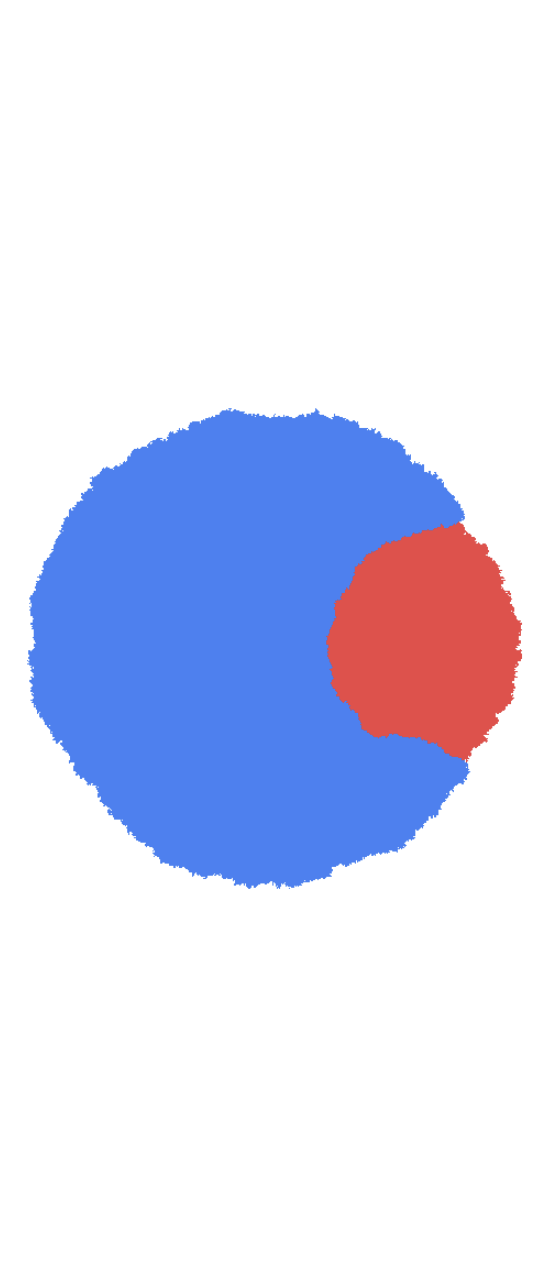}\hspace{2mm}%
        \includegraphics[height=3cm]{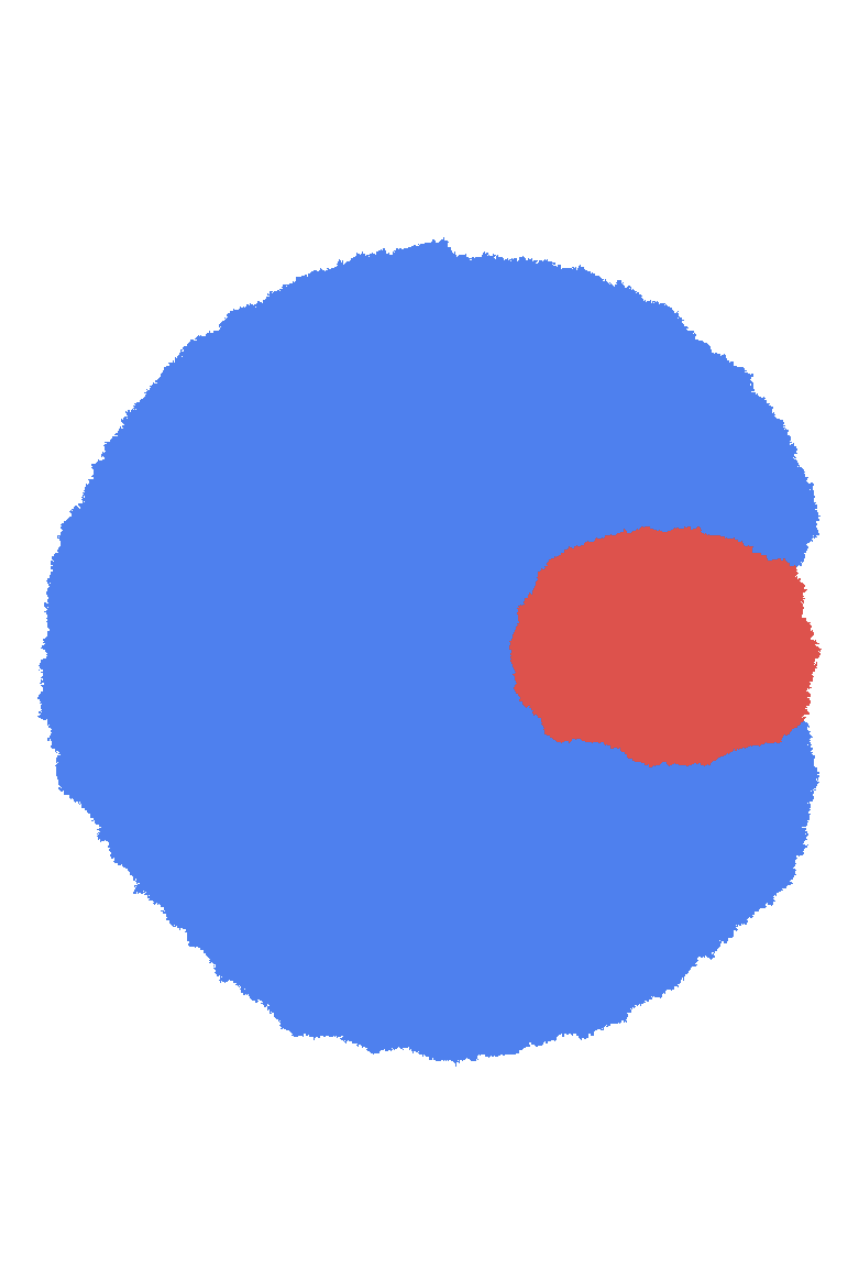}\hspace{2mm}%
        \includegraphics[height=3cm]{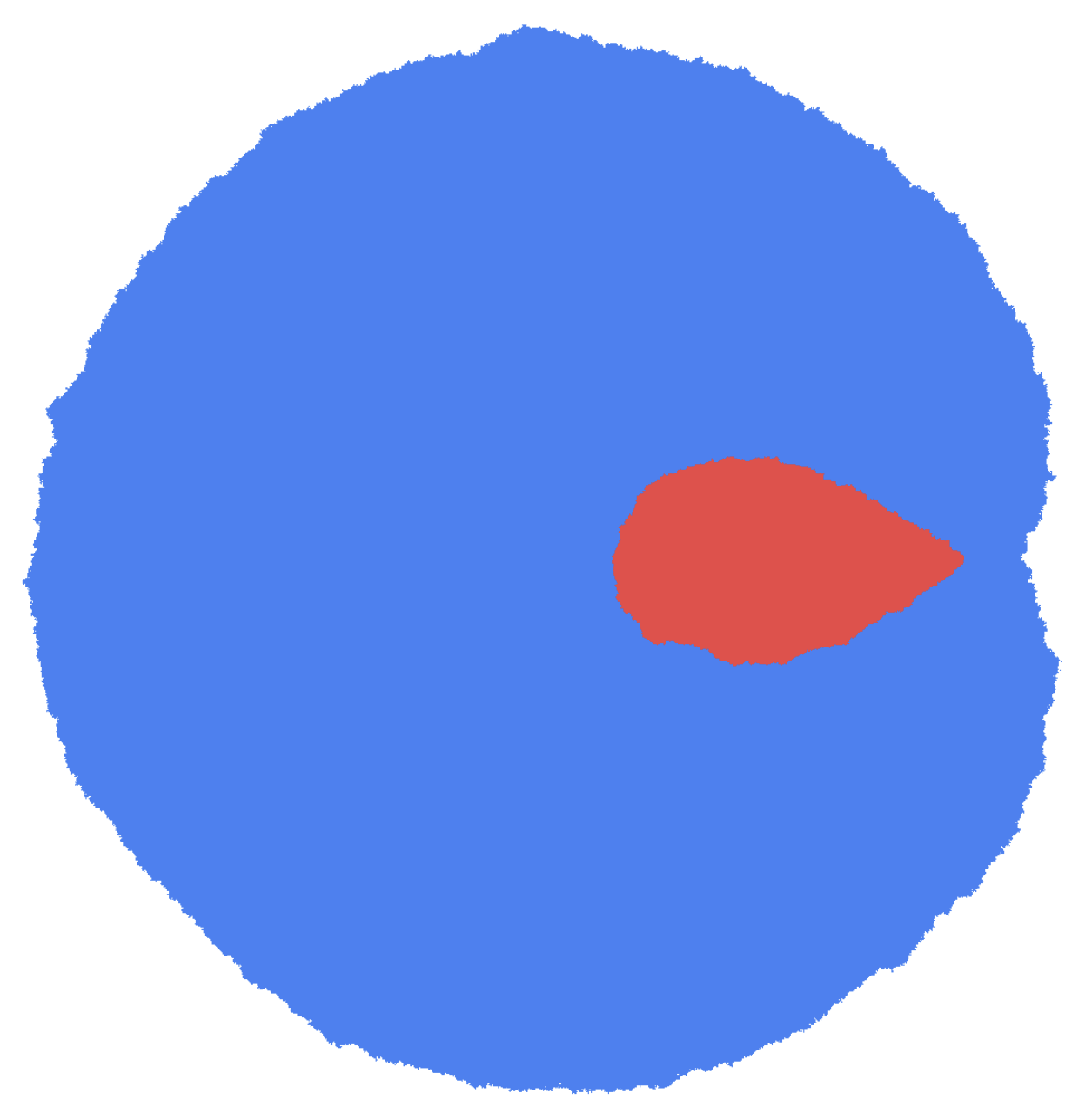}%
    }
    \caption{A simulation of the two-type Richardson model with rates $1$ and $1.5$. Snapshots taken at times (left to right) 20, 50, 100, 170 and 250.}
    \label{fig:five_images}
\end{figure}

\section{Introduction}

We study the two-type Richardson model introduced by H{\"a}ggstr{\"o}m and Pemantle~\cite{haggstrom1998first}. In this model two competing species grow on the integer lattice $\mathbb Z^2$. The process can be viewed as a continuous time Markov chain $\{(W_t,S_t)\}_{t\ge 0}$ where $W_t,S_t\subseteq \mathbb Z^2$ are disjoint, growing subsets of $\mathbb Z^2$ starting from $W_0=\{0\}$ and $S_0=\{v_0\}$ with $v_0\in \mathbb Z^2\setminus \{0\}$. Given two infection rates $\lambda _1,\lambda _2>0$,  the vertices in $W_t$ infect healthy vertices (i.e., vertices outside of $ W_t\cup S_t$) at rate $\lambda _1$ and the vertices in $S_t$ infect healthy vertices at rate $\lambda _2$. With positive probability one type will encircle the other and make it stop growing (see Figure~\ref{fig:five_images}). The question is whether coexistence can occur in which both species continue to grow forever. We denote this event by $\mathscr C:=\{|W_t|\to \infty $ and $|S _t|\to \infty \}$. A longstanding conjecture due to H{\"a}ggstr{\"o}m and Pemantle~\cite{haggstrom1998first} states that $\mathbb P(\mathscr C)>0$ if and only if $\lambda _1=\lambda_2$. See the surveys \cite{deijfen2015pleasures,blair2010first}  and \cite[Section~6]{auffinger201550} about this conjecture. While the ``if'' part of the conjecture is known \cite{haggstrom1998first,garet2005coexistence,hoffman2005coexistence,deijfen2007two}, there are only partial results toward the ``only if'' part \cite{haggstrom2000absence,garet2008first,ahlberg2020two}. See Section~\ref{sec:previous} for more details. In this paper we resolve this conjecture.

\begin{thm}\label{thm:main}
    For any $\lambda_1 \neq \lambda _2$ we have that $\mathbb P (\mathscr{C})=0$.
\end{thm}

\begin{remark}
    The proof in this paper uses the recent (unverified) result by OpenAI that the limit shape of exponential FPP is strictly convex (see Theorem~\ref{thm:openai} below). This result made the proof much simpler, but it is not crucial for our argument. In fact, before the announcement of OpenAI on October 6, 2026, we were writing a first version of this paper in which we deal with possible flat edges of the limit shape. In a future work (or a future version of this paper), we intend to pursue our old proof in order to generalize the result to competition models with general weight distributions, for which strict convexity is still open.
\end{remark}

\subsection{Previous works}\label{sec:previous}

H{\"a}ggstr{\"o}m and Pemantle~\cite{haggstrom1998first} introduced the model on $\mathbb Z^2$ in 1997 and made the conjecture that coexistence can occur if and only if $\lambda _1=\lambda _2$. In the same paper they proved the if part of the conjecture. This was later extended to higher dimensions and to other weight distributions by Garet and Marchand~\cite{garet2005coexistence} and independently by Hoffman~\cite{hoffman2005coexistence}. An additional proof of this result was given by Deijfen
and H{\"a}ggstr{\"o}m~\cite{deijfen2007two}.
The same authors~\cite{deijfen2006initial} proved that on $\mathbb Z^d$ the answer to the question of whether coexistence occurs with positive probability or not, is independent of the initial configuration $W_0,S_0$. 

Significant progress towards the only if part of the conjecture has been made by H{\"a}ggstr{\"o}m and Pemantle~\cite{haggstrom2000absence} who proved that in any dimension coexistence occurs with zero probability for all ratios of $\lambda _1/\lambda _2$ except for at most countably many such ratios. This was extended to general weights by Garet and Marchand~\cite{garet2008first} who also proved that if coexistence happens when the rates are not equal, then only a negligible fraction of space is occupied by the strong type. Ahlberg, Deijfen, and Hoffman~\cite{ahlberg2020two} proved that the only if part of the conjecture holds in the half space $\mathbb Z\times \mathbb Z_+$ when the initial infections start from two points on the boundary. In this setup they showed that if the strong type survives then it necessarily occupies a ray along the boundary on which it moves faster and eventually circles the weak type.

We refer the reader to work on unbounded initial configurations
\cite{deijfen2007two,blair2012geometric,ahlberg2026surviving},
finitely many competing types
\cite{hoffman2008geodesics,gouere2007shape,pimentel2007multitype},
and infinitely many types
\cite{damron2013examples,auffinger2013differentiability,antunovic2017stationary,ahlberg2021existence}.
See also~\cite{deijfen2004stochastic,deijfen2004coexistence,gouere2007shape,carstens2011two} for continuous space analogues
 of the model and 
\cite{garet2005coexistence,deijfen2006nonmonotonic,pimentel2007multitype,ahlberg2020two,candellero2021coexistence} for the analysis of the model on infinite graphs beyond $\mathbb Z^d$. Finally, we refer the reader to \cite{benjamini2002survival,MR4483139} for a deterministic version of the model and to~\cite{garet2006competition,sidoravicius2019multi,candellero2021coexistence,finn2024nonequilibrium,candellero2024hostile,finn2022conversion} for additional closely related models.

\subsection{Conventions}

By scaling of time and possibly switching the roles of $W$ and $S$, we may and will assume throughout the paper that
\begin{equation}
    \lambda _1=1 \quad \text{ and }\quad \lambda _2=\lambda >1.
\end{equation}
So $W$ stands for the Weak type and $S$ stands for the Strong type. We will always refer to the strong and weak types as blue and red respectively. These colors will also appear in all the figures.

We use $C,c$ to denote positive constants depending only on $v_0$ and $\lambda $. These constants are regarded as generic in the sense that their value may change from one appearance
to the next.

\subsection{Sketch of proof and key propositions}

The first relevant observation is a topological one. By planarity and using the fact that both $W_t$ and $S_t$ are connected, it follows that the boundary of the infected region $W_t\cup S_t$ is the union of two disjoint arcs, the blue part of the boundary and the red part. This makes the geometric picture simpler in two dimensions and is the main difficulty to overcome when $d\ge 3$ (see Section~\ref{sec:open}).

\begin{figure}
    \centering
        \begin{overpic}[width=0.6\linewidth]{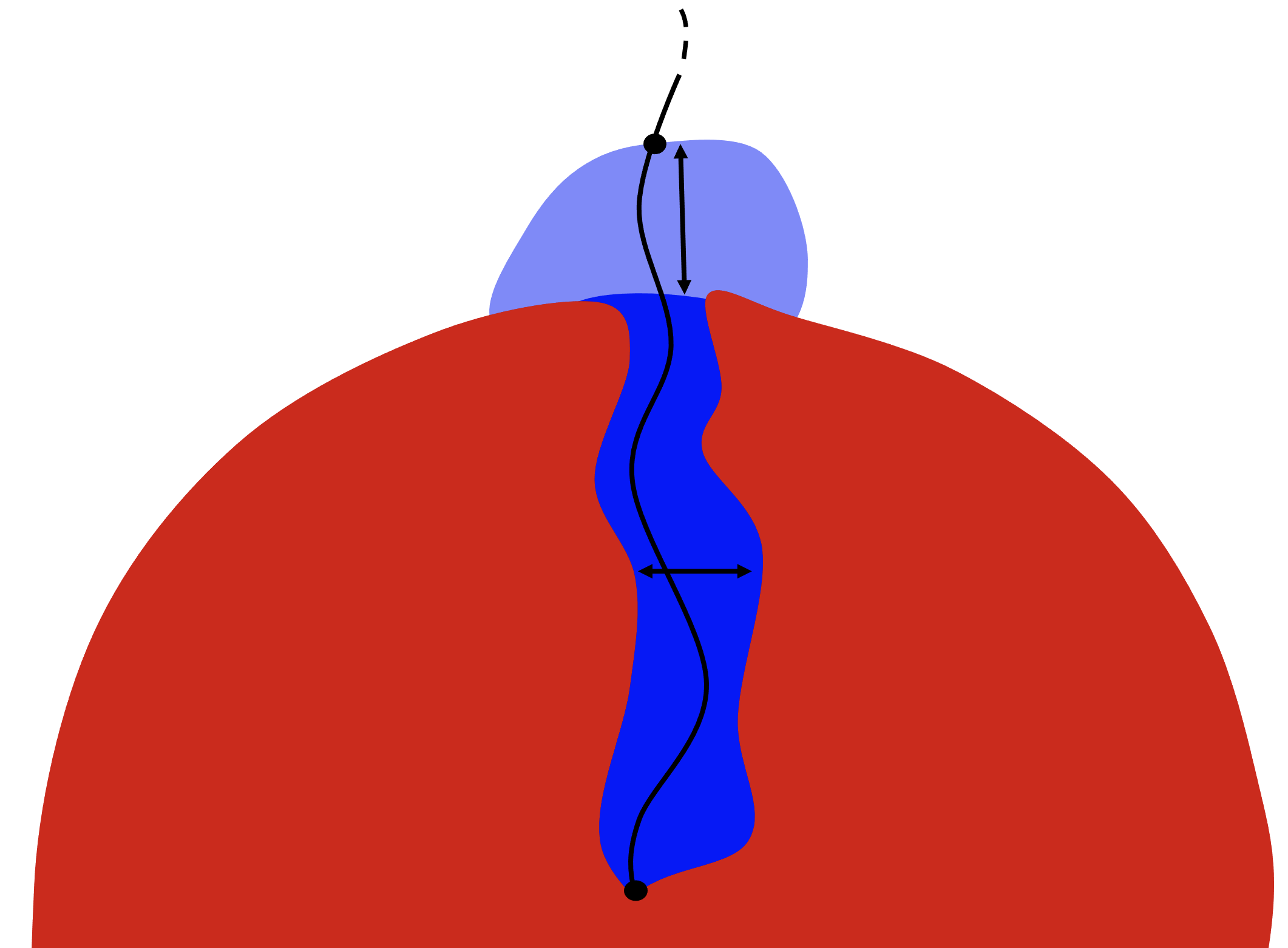}
        \put(46,2.1){\ensuremath{v_0}}
        \put(54,55){\ensuremath{T^{\frac{3}{8}}}}
        \put(46.5,64){\ensuremath{u_k}}
        \put(52,30){\ensuremath{T^{\frac{1}{4}}}}
    \end{overpic}
    \caption{Proof of Lemma~\ref{lem:inclusion}.}
    \label{fig:mermin}
\end{figure}

Next, it is not hard to check that if blue is very far outside the convex hull of red, then with high probability it will move faster, escape red and eventually encircle red. Hence, the most difficult scenario to rule out is when red almost circles blue completely, and only leaves a narrow and winding tunnel for blue to escape through (see, e.g., Figure~\ref{fig:mermin}). Even though blue is faster, the winding of the tunnel left by red will make it travel a longer distance. If this is done exactly right, by the time blue exits this tunnel, red will already create a larger tunnel for blue, and this can potentially continue forever in all scales. It sounds completely absurd that coexistence will occur in this way since red is growing randomly and is not creating tunnels for blue adversarially. However, it is not clear how to prove this rigorously, and in fact our proof uses the random fluctuations of blue and not of red to rule this out.

In the proof we establish two key propositions. Take a large time $T$ and consider the size of the blue part of the boundary up to time $T$. This is the width of the blue tunnel. The propositions say that:
\begin{enumerate}
    \item 
    If the size of this boundary is always $\lesssim T^{1/4}$, then with high probability red will eventually encircle blue completely.
    \item If at some point, the size of this boundary is $\gg T^{1/4}$, then with high probability blue will escape and encircle red.
\end{enumerate}

Let us state these key propositions precisely. Fix $\epsilon :=0.01$ throughout this paper. Let $\mathcal L_t$ be the blue part of the boundary at time $t$ (not including the part of the boundary in holes):
\begin{equation}
    \mathcal L_t := \big\{ u\in S_t : \exists v\sim u \text{ in the infinite connected component of }  \mathbb Z^2 \setminus (W_t\cup S_t) \big\}.
\end{equation}
Let $L_t$ be the diameter of $\mathcal L_t$:
\begin{equation}
    L_t:=\max \big\{ \|u-v\| : u,v\in \mathcal L_t  \big\}.
\end{equation}

It is easier to work with the diameter of the boundary rather than the size of the boundary because the diameter is more smooth in time. We write $W_\infty =\bigcup_{t\ge 0} W_t$ and $S_\infty =\bigcup_{t\ge 0} S_t$.

\begin{prop}\label{prop:mermin}
    We have that
    \begin{equation}
        \mathbb P \Big( |S_\infty |=\infty   \ \ \text{and } \   \forall t\le T, \  L_t \le T^{1/4+3\epsilon } \Big) \to 0,\quad T\to \infty .
    \end{equation}
\end{prop}

\begin{prop}\label{prop:limit}
    We have that
    \begin{equation}
        \mathbb P \Big(  |W_\infty |=\infty   \ \ \text{and } \exists t\le T,\  L_t \ge T^{1/4+3\epsilon }  \Big) \to 0,\quad T\to \infty  .
    \end{equation}
\end{prop}

Theorem~\ref{thm:main} easily follows from these propositions. Indeed, combining the two propositions  and taking the limit $T\to \infty $ we obtain that $\mathbb P (\mathscr C)=\mathbb P (|W_\infty |=\infty \text{ and }|S_\infty |=\infty )=0$.

Let us briefly explain the ideas in the proof of these propositions. Proposition~\ref{prop:mermin} is proved in Section~\ref{sec:mermin} using a Mermin--Wagner style argument. Such arguments were used before in first and last passage percolation \cite{dembin2024coalescence,dembin2025influence,elboim2025small,dembin2025noise,schmid2026mixing,chatterjee2019general,bates2020fluctuation} but not in an adaptive way as in here. The idea is very simple. Since the blue tunnel is very narrow, the blue speed along this tunnel depends on a small number of random variables. The random fluctuations of these variables can push the blue inside and kill it, or push it outside and help it escape. More precisely, there are $n\approx T^{5/4}$ exponential variables in this tunnel (length $T$ and width $T^{1/4}$) that determine the speed of blue (see Figure~\ref{fig:mermin}). With non-negligible chance all these variables will be smaller by $n^{-1/2}\approx T^{-5/8}$ (see Lemma~\ref{lem:mermin}) pushing blue distance $\approx T^{3/8}$ forward along its tunnel and making its boundary much wider than $T^{1/4}$. One small difficulty is that this tunnel is random and we have to choose which variables to decrease in an evolving and predictable way. 

We turn to describe the proof of Proposition~\ref{prop:limit} which is done in Section~\ref{sec:limit}. This part uses the connection of the two-type Richardson model with exponential first-passage percolation (FPP). We analyze the geometry of the FPP ball $B_t$ at time $t$, which is closely related to the infected region in our process. See Section~\ref{sec:FPP} for more details about this connection. The idea is to consider a second order approximation $\tilde B_t$ of the FPP ball $B_t$. See \eqref{eq:second} and Figure~\ref{fig:second}. A first order approximation of this ball is simply $t$ times the limit shape, with an error of order $\sqrt{t}$, which is too large for our purposes. In the second order approximation $\tilde B_t$ we consider the random FPP ball at time $t-\sqrt{t}$ and ``inflate'' each of its points by $\sqrt{t}$ times the limit shape. This approximation is much closer to the real ball $B_t$ with an error of order $t^{1/4}$. Nevertheless, this approximation is more regular than the real ball. We show that all the points on the boundary of $\tilde B_t$ are sufficiently exposed so that a blue infection starting from the boundary of $\tilde B_t$ will manage to escape a red infection that is contained in $\tilde B_t$. This part uses the strict convexity of the limit shape that was recently uploaded by OpenAI. Assuming that the blue boundary is of width $\gg t^{1/4}$, this will be enough to overcome the random fluctuations of order $t^{1/4}$ around the second order approximation, allowing blue to escape the real infected region.

\begin{figure}[htbp]
    \centering
    \includegraphics[width=0.45\linewidth]{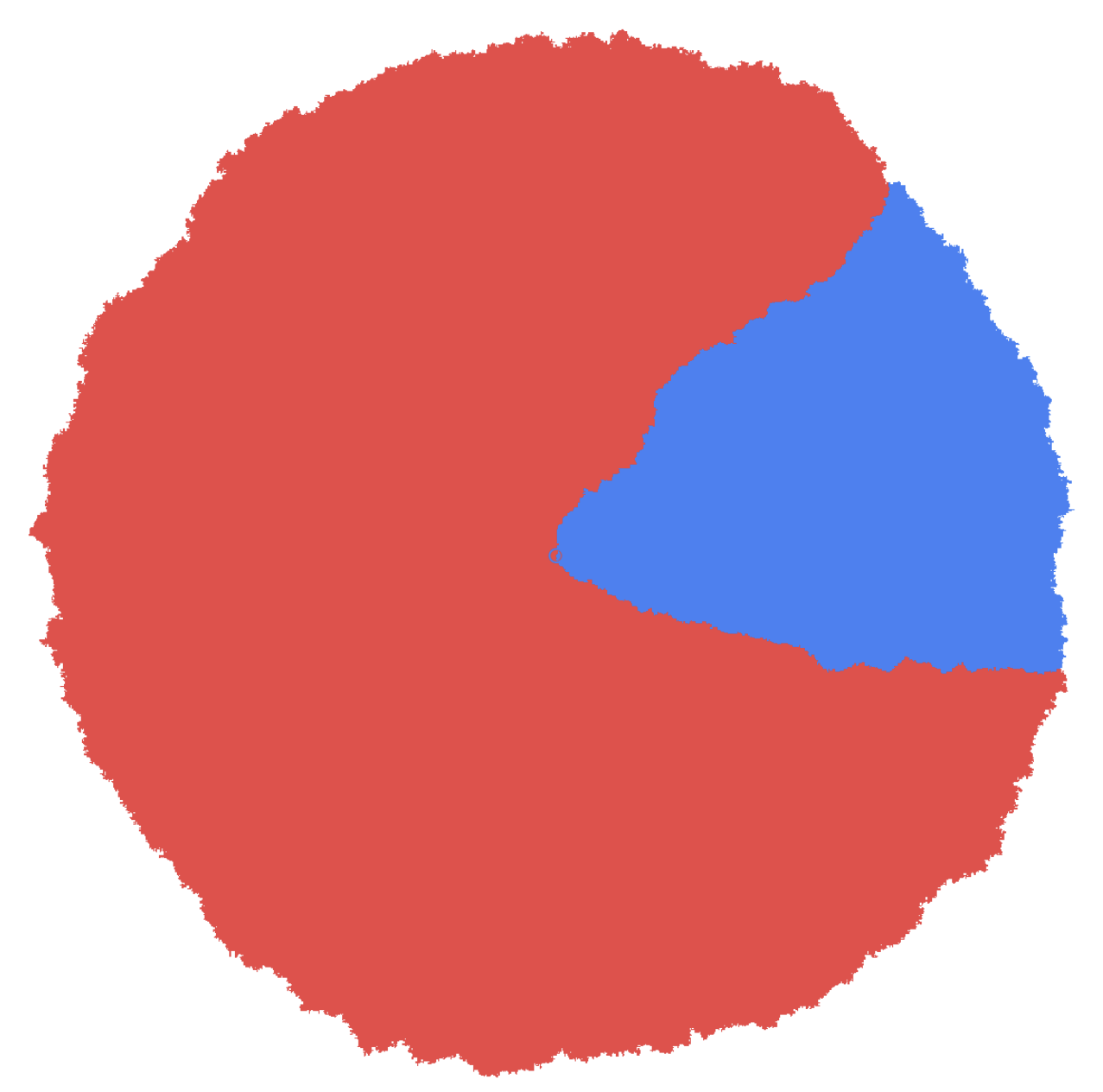}
    \hfill
    \includegraphics[width=0.45\linewidth]{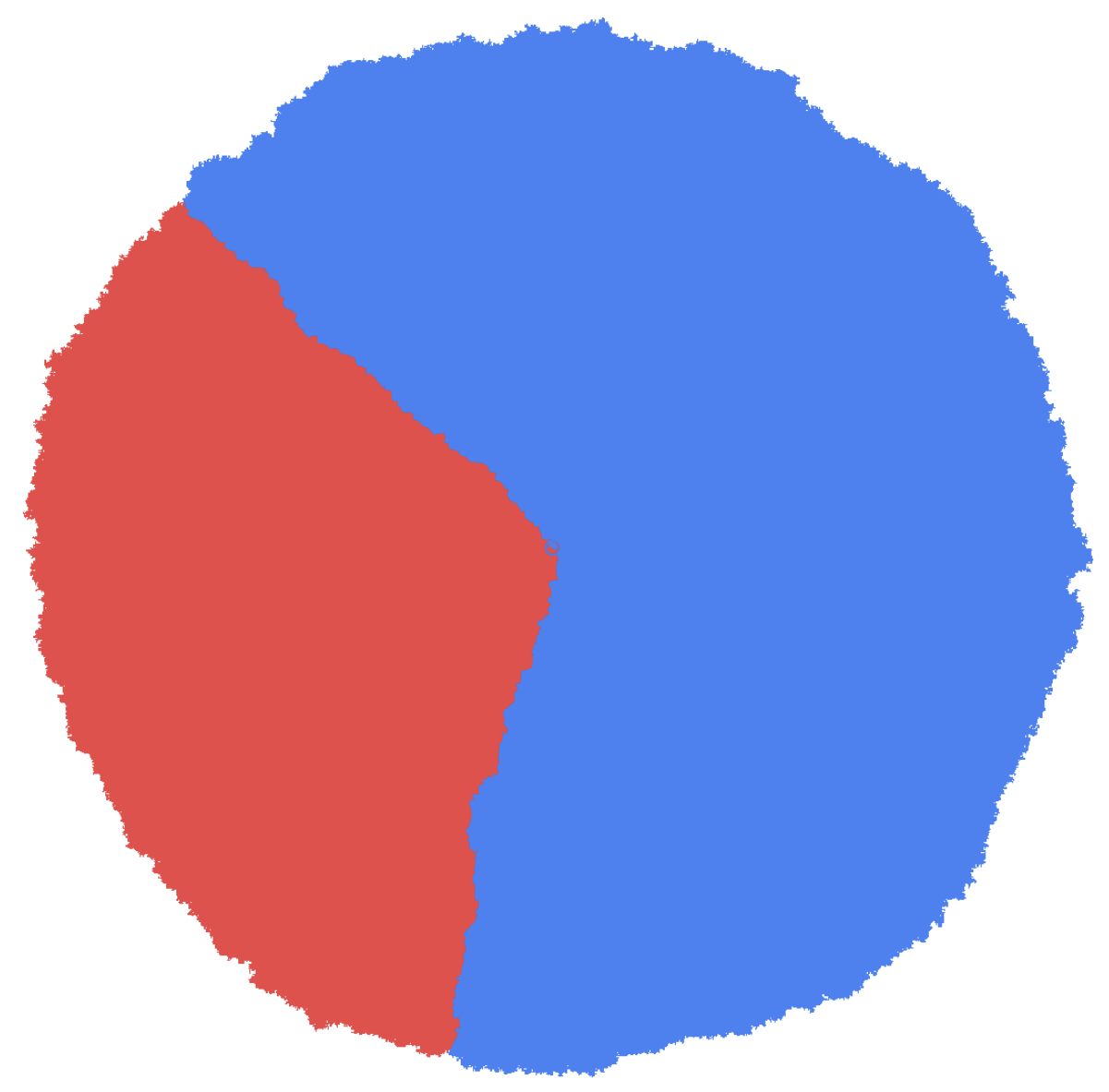}
    \caption{A simulation of the two-type Richardson model with equal rates.}
    \label{fig:equal}
\end{figure}

\subsection{Open Problems}\label{sec:open}

The most interesting open problem is to extend the result to all dimensions $d\ge 3$. In dimension $d\ge 3$ there is no single blue tunnel. Blue can have several ``tentacles'' connecting $v_0$ to the boundary of $W_t\cup S_t$. As time varies, these tentacles can split and merge, and blue can decide from which tentacle to escape only in the future. This makes the geometric picture much more complicated, and it is not clear how to implement our Mermin--Wagner approach discussed above.

Another interesting question concerns the two-type Richardson model with equal rates. In this case coexistence occurs with positive probability, and on the event of coexistence we would like to understand what the sets $W_t$ and $S_t$ look like for large $t$. Their union $W_t\cup S_t$ divided by $t$ will converge to the limit shape of exponential FPP. It is natural to expect that $W_t/t$ and $S_t/t$ will converge almost surely to a random ``pizza slice'' and its complement (see Figure~\ref{fig:equal}). Such results were obtained for a rotationally symmetric version of this model \cite{pimentel2007multitype}. It is very likely that one can prove it on the lattice using the new result of OpenAI on the strict convexity of the limit shape of exponential FPP.

\subsection{Acknowledgments and AI statement}

We thank Ofer Busani, Ilay Hoshen, Gur Lifshitz and Michael Chapman for fruitful discussions.

All of the main ideas in this paper were generated by humans. We used ChatGPT to write the simulations in Figure~\ref{fig:five_images} and Figure~\ref{fig:equal} and to review our paper. It also simplified our proofs in Lemma~\ref{lem:mermin}, Lemma~\ref{lem:limitshape} and the construction of $\gamma _1,\gamma _2$ in Lemma~\ref{lem:kill}.

\section{preliminaries}\label{sec:FPP}

The one-type and the two-type Richardson models are closely related to first-passage percolation with exponential weights. Let us review some basic facts from FPP.

\subsection{Exponential FPP}

Let $E(\mathbb Z ^2)$ be the set of nearest neighbor edges of $\mathbb Z^2$ and let $\{X_e\}_{e\in E(\mathbb Z^2)}$ be i.i.d.\ $\exp(1)$ edge weights. For two vertices $u,v\in \mathbb Z^2$ define the passage time between them to be
\begin{equation}
    T(u,v):=\min _{p} \sum _{e\in p} X_e, 
\end{equation}
where the minimum is over $\mathbb Z^2$ paths $p$ connecting $u$ to $v$. This gives a random metric on $\mathbb Z^2$. It is sometimes convenient to extend this metric  to $\mathbb R^2$ by letting $T(x,y)=T(u,v)$ for every $x\in u+[-\tfrac{1}{2},\tfrac{1}{2})^2$ and $y\in v+[-\tfrac{1}{2},\tfrac{1}{2})^2$. We let 
\begin{equation}
    B_t:=\{x\in \mathbb R^2 : T(0,x) \le t\}.
\end{equation}

The Cox--Durrett limit shape theorem~\cite{cox1981some} states that there is a deterministic origin symmetric convex shape $\mathcal B \subseteq \mathbb R^2$ such that the ball $B_t$ satisfies $B_t/t\to \mathcal B$ (in the Hausdorff distance). So that $t\mathcal B$ gives a ``first order approximation'' of the ball $B_t$. In our proof we will also consider a better, second order approximation of $B_t$.

\begin{thm}[OpenAI]\label{thm:openai}
    The limit shape $\mathcal B$ is strictly convex.
\end{thm}

We remark that if the distribution of the edge weights is not exponential, strict convexity does not necessarily hold. See, e.g., \cite{flat} and the related \cite{damron2013examples}.

Let $\mu :\mathbb R ^2\to [0,\infty )$ be the limiting norm of exponential FPP, which is the norm whose unit ball $\{x\in \mathbb R^2: \mu (x)\le 1\}$ is equal to $\mathcal B$. For two vertices $u,v\in \mathbb Z^2$ the norm $\mu (u-v)$ gives a first order approximation for the passage time $T(u,v)$.

\begin{remark}\label{rem:1}
    Strict convexity and a compactness argument imply that for any $\delta >0$, there exists $\eta >0$ depending only on $\mathcal B$ and $\delta $ such that for any two directions $x,y\in \partial \mathcal B$ with $\mu \big(\frac{x+y}{2}\big)\ge 1-\eta $ we have that $\mu (x-y)\le \delta $.
\end{remark}

We turn to discuss the fluctuations of $T(u,v)$ around $\mu (u-v)$. Alexander's bound on the non-random fluctuations~\cite[Theorem~3.2]{alexander1997approximation} states that
\begin{equation}\label{non-random}
    \big|\mathbb E [T(u,v)]-\mu (u-v) \big| \le C\sqrt{\|u-v\|} \log (\|u-v\|).
\end{equation}
Moreover, it follows from, e.g., \cite[Theorem 3.11]{auffinger201550} that
\begin{equation}\label{random}
    \mathbb P \Big(  \big|T(u,v)-\mathbb E [T(u,v)] \big| \ge r \sqrt{\|u-v\|} \Big) \le Ce^{-cr}. 
\end{equation}
In our proof at time scale $T$, it is convenient to work on the event $\mathcal E =\mathcal E _T$ defined by 
\begin{equation}\label{eq:E}
    \mathcal E := \Big\{ \forall u,v \in T^2\mathcal B \cap \mathbb Z^2, \ \big|T(u,v)-\mu (u-v) \big| \le T^{\epsilon /9} \sqrt{\|u-v\|}  \Big\}.
\end{equation}
The following lemma is an immediate consequence of \eqref{random} and \eqref{non-random} and a union bound.

\begin{lem}\label{lem:E}
    We have that $\mathbb P (\mathcal E)\ge 1-e^{-cT^c}$.
\end{lem}

\subsection{The natural coupling with exponential FPP}\label{sec:natural}

We can sample our process using the edge weights $\{X_e\}_{e\in E(\mathbb Z^2)}$. Indeed every time $W_t\cup S_t$ reaches an endpoint of an edge $e$ for the first time, it uses $X_e$ to sample the infection time through $e$. If red reaches the endpoint, it infects the other endpoint after time $X_e$, and if blue reached the endpoint, it infects the other endpoint after time $X_e/\lambda \sim \exp (\lambda )$. In this construction every edge weight is examined by the process only once and so different infection times along edges are independent. Hence this construction recovers the two-type Richardson process.

\begin{lem}
    In the natural coupling with exponential FPP we have that
    \begin{equation}\label{eq:WBWS}
    W_t\subseteq B_t\subseteq W_t\cup S_t,
    \end{equation}
\end{lem}
(One may understand $W_t,S_t,B_t$ either as subsets of $\ZZ^2$ or as continuous sets by adding $[-\frac{1}{2},\frac{1}{2})^2$ around each integer point.)

\begin{proof}
    It suffices to prove that $T^-(v)\le T(v)\le T^+(v)$ for any $v\in \mathbb Z^2$, where 
    \begin{equation*}
        T^-(v):=\inf \{t:v\in W_t\cup S_t\},\quad T(v):=\inf \{t: v\in B_t\},\quad T^+(v):=\inf \{t:v\in W_t\}.
    \end{equation*}
    To prove the first inequality (equivalently, the second containment) observe that for any edge $e=\{u,v\}$ we have that $T^-(v)\le T^-(u) + X_e$,
    since once $u$ is infected, it will infect $v$ within time $X_e$ if $u\in W_t$ and within time $X_e/\lambda <X_e$ if $u\in S_t$. This gives
\begin{equation}
    T^-(v)\le \inf _{p} \sum _{e\in p} X_e =T(v),
\end{equation}
where the infimum is over paths connecting $0$ to $v$.

To prove the second inequality, let $p$ be the path of infections from $0$ to $v$ of the weak species. We have that 
\begin{equation*}
    T(v) \le \sum _{e\in p} X_e =T^+(v).\qedhere
\end{equation*}
\end{proof}

\section{Mermin--Wagner: proof of Proposition~\ref{prop:mermin}}\label{sec:mermin}

In this section we prove Proposition~\ref{prop:mermin}. We would like to prove that $\mathbb P (F\cap G)\to 0$ where
\begin{equation}
    F:=\{|S_{\infty }|=\infty \}  \quad \text{ and } \quad  G:= \big\{ \forall t\le T, \  L_t \le T^{1/4+3\epsilon }  \big\}.
\end{equation}
This is done using a Mermin--Wagner type argument. The idea is to decrease slightly the infection times of the strong type, making it a bit faster, and making its boundary larger. To keep the distribution of the perturbed process sufficiently close to the original one, we only accelerate the infection times along edges when the boundary $\mathcal{L}_t$ is small. 

Recall the natural coupling from Section~\ref{sec:natural} in which the competition process $P$ is sampled using the exponential FPP environment $\{X_e\}_{e\in E(\mathbb Z^2)}$. Next, using this environment we define the perturbed process $P'=\{(W_t',S_t')\}_{t\ge 0}$  that is coupled with $P$. In this process we will use prime to denote the corresponding objects such as $\mathcal L_t'$ or $L_t'$. 

The process $P'$ and the perturbed environment $\{X_e'\}_{e\in E(\mathbb Z^2)}$ are defined inductively. Suppose that $P'$ was defined up to time $t$ and $t$ is the first time an endpoint $u$ of an edge $e$ is infected in $P'$. We define $X_e'$ by speeding slightly the blue infection whenever its boundary is not too wide. Formally,
\begin{equation}
\delta :=T^{-5/8-3\epsilon } \quad \text{ and }\quad      X_e':= \begin{cases}
        \frac{X_e}{1+\delta } , \quad \text{if } t\le T, \ u\in \mathcal L_t', \ L_t'\le T^{1/4+ 3\epsilon } \\
        X_e, \quad \text{otherwise}
    \end{cases}.
\end{equation}
Now, as in the original process, the infection time along $e$ in $P'$ will be $X_e'$ if $u\in W_t'$ and $X_e'/\lambda $ if $u\in S_t'$.

Proposition~\ref{prop:mermin} will follow from the next two lemmas. To state the first lemma, recall the definition of $\mathcal E$ from \eqref{eq:E} and let $F '$ and $G'$ be the events corresponding to $F$ and $G$ in the perturbed process $P'$. The lemma says that if blue survives, then after making it faster, its boundary will be wide.

\begin{lem}\label{lem:inclusion}
    We have that $F \cap \mathcal E \subseteq F' \setminus G '$, if $T$ is sufficiently large.
\end{lem}

The second lemma says that the perturbation did not change the distribution of the competition process $P$ by much.

\begin{lem}\label{lem:TV}
    We have that $d_{\rm TV}(P,P')\le CT^{-\epsilon }$.
\end{lem}

Using these two lemmas, we can easily prove Proposition~\ref{prop:mermin}.

\begin{proof}[Proof of Proposition~\ref{prop:mermin}]
    By Lemma~\ref{lem:inclusion} we have that 
\begin{equation}
\begin{split}
        \mathbb P (F) -\mathbb P (\mathcal E^c) &\le \mathbb P ( F \cap \mathcal E) \le  \mathbb P ( F' \setminus G ')\\
        &\le \mathbb P ( F \setminus G ) +d_{\rm TV} (P,P')=\mathbb P ( F )-\mathbb P ( F \cap G ) +d_{\rm TV} (P,P').
\end{split}
\end{equation}
Rearranging and using that  $\mathbb P (\mathcal E^c)\to 0$ by Lemma~\ref{lem:E} and that $d_{\rm TV} (P,P') \to 0$ by Lemma~\ref{lem:TV} we obtain $\mathbb P (F \cap G )\to 0$.
\end{proof}

The rest of this section is devoted to the proof of Lemma~\ref{lem:inclusion} and Lemma~\ref{lem:TV}.

\subsection{The effect of the perturbation: Proof of Lemma~\ref{lem:inclusion}}

For the proof of Lemma~\ref{lem:inclusion} we will need the following claim.

\begin{claim}\label{claim:monotonicity}
    For all $t\ge 0$ we have that $S_t\subseteq S_t'$ and $W_t\supseteq W_t'$.
\end{claim}

Claim~\ref{claim:monotonicity} can be proved using a straightforward induction on time. The details are omitted.

\begin{proof}[Proof of Lemma~\ref{lem:inclusion}]
It follows from Claim~\ref{claim:monotonicity} that $F \subseteq F'$ and therefore it suffices to prove that on $F\cap \mathcal E$ the event $G'$ does not hold. Suppose by contradiction that $G'$ holds.

Let $p=(v_0,v_1,\dots )$ be an infinite path of blue 
infections in the original environment. Such a path exists on $F$. Let $T(v_j)$ be the time blue infected $v_j$. Let $T'(v_j)$ be the time blue infects $v_j$ in the perturbed environment. By Claim~\ref{claim:monotonicity} we have $T(v_j)\ge T'(v_j)$.

We define a sequence of vertices $u_0,u_1,\dots $ inductively as follows (it will be a subsequence of $v_0,v_1,\dots $). Let $u_0=v_0$ and for $i\ge 0$ define $u_{i+1}$ to be the first vertex to be infected in the perturbed process among all vertices lying strictly after $u_i$ along $p$. Let $u_k$ be the last vertex along this sequence with $T'(u_k)\le T$. It suffices to prove that on $F\cap \mathcal E\cap G'$ we have
\begin{equation}\label{eq:far}
    d_{\mu }(u_k, W'_{T})\ge  \tfrac{1}{6}T^{3/8-3\epsilon }.
\end{equation}
Indeed, by construction at time $T$ all vertices after $u_k$ along $p$ are healthy and so $u_k\in \mathcal L _{T}'$. Hence \eqref{eq:far} implies that $L'_{T}\ge cT^{3/8-3\epsilon }>T^{1/4+3\epsilon }$ which contradicts the fact that $G'$ holds.

We turn to prove \eqref{eq:far}. First, note that $T-T'(u_k) \le T^{\epsilon }$. Indeed, once $u_k$ is infected, it will infect the next vertex $u_k^+$ along $p$ in time at most $T'(u_k^+)-T'(u_k) \le T(u_k,u_k^+)\le T^{\epsilon }$ on the event $\mathcal E$, and by definition of $u_k$ we have that $T'(u_k^+)>T$.
Next, note that if some $u_i$ is infected at time $t\le T$ then $u_i\in \mathcal L_t'$. It follows that the infection time along the edge $e$ connecting $u_i$ to the next vertex on $p$ is $X_e/(\lambda (1+\delta ))$. Thus, for all $i\le k$ we have that  $T'(u_{i+1})-T'(u_i) \le (T(u_{i+1})-T(u_i))/(1+\delta )$.  This gives that
\begin{equation}
    T -T^{\epsilon } \le T'(u_k) \le \frac{T(u_k)}{1+\delta }  \le  T(u_k)- \tfrac{1}{4}T^{3/8-3\epsilon },
\end{equation}
where in the last inequality we used that $\frac{1}{1+x}\le 1-x/2$ for $x\in [0,1]$ and that $T(u_k)\ge T'(u_k)\ge T-T^{\epsilon }$ on $\mathcal E$. Rearranging, we get that $T(u_k)\ge T+\frac{1}{5}T^{3/8-3\epsilon }$. However, if it takes another $\frac{1}{5}T^{3/8-3\epsilon }$ time for $S_{T}\cup W_{T}$ to reach $u_k$, then it means that on $\mathcal E$ we have $d_{\mu }(u_k, S_{T}\cup W_{T}) \ge \frac{1}{6}T^{3/8-3\epsilon }$. Since $W_{T}'\subseteq W_{T}$, this finishes the proof of \eqref{eq:far}, and hence also of Lemma \ref{lem:inclusion}.
\end{proof}

\subsection{The total variation distance: Proof of Lemma~\ref{lem:TV}}

\begin{lem}\label{lem:mermin}
    Let $n\ge 1$ be an integer and let $\delta >0$. Let $Z_1,\dots ,Z_n$ be i.i.d.\ $\exp (1)$ random variables. For $i\le n$ let $Z_i':=Z_i/(1+\delta )\sim \exp (1+\delta )$. Let $Z=(Z_1,\dots ,Z_n)$ and $Z'=(Z_1',\dots ,Z_n')$. Then
    \begin{equation}
        d_{{\rm TV}}(Z,Z') \le \delta \sqrt{n}.
    \end{equation}
\end{lem}

\begin{proof}
For $t\ge 1$, define
\[
f_t(x)=t^n e^{-t\sum_{i=1}^n x_i},
\qquad x\in\mathbb R_+^n.
\]
The densities of $Z$ and $Z'$ are $f_1$ and $f_{1+\delta }$ respectively. We have that
\[
\frac{\partial f_t(x)}{\partial t}
=\left(\frac nt-\sum _{i=1}^nx_i\right)f_t(x)
\]
and therefore letting $S$ be the sum of $n$ i.i.d.\ $\exp(t)$ random variables, we have that
\[
\int_{\mathbb R_+^n}
\left|\frac{\partial f_t(x)}{\partial t}\right|\,dx
=\mathbb E \left|S-\frac nt\right|
\le \sqrt{\operatorname{Var}(S)}
\le \sqrt{n},
\]
where here we used that $\mathbb E [S]=n/t$ and the Cauchy-Schwarz inequality. We obtain
\begin{align*}
2d_{{\rm TV}}(Z,Z')
&=\int_{\mathbb R_+^n}
  |f_{1+\delta }(x)-f_1(x)|\,dx\\
&=\int_{\mathbb R_+^n}
  \left|\int_1^{1+\delta }
  \frac{\partial f_t(x)}{\partial t}\,dt\right|\,dx \le \int_1^{1+\delta }
  \int_{\mathbb R_+^n}
  \left|\frac{\partial f_t(x)}{\partial t}\right|\,dx\,dt\le  \delta \sqrt n.\qedhere
\end{align*}
\end{proof}

Let $\{X_e\}_{e\in E(\mathbb Z^2)}$ be the environment of exponentials. Let $P=\{(W_s,S_s)\}_{s\ge 0}$ be the original process driven by the environment $X$.

For the proof of Lemma~\ref{lem:TV} we will need the following claim. To state the claim, let $N':=|\{e:X_e'< X_e\}|$  be the 
number of edge weights that were decreased in $P'$.

\begin{claim}\label{claim:k}
    We have that $\mathbb P \big( N' \ge CT^{5/4+3\epsilon } \big) \le e^{-cT^c}$.
\end{claim}

\begin{proof}
    By Lemma~\ref{lem:E} it suffices to prove that on $\mathcal E$ we have that $N' < CT^{5/4+3\epsilon }$.  If for some $t\le T$, we have that $L_t'\le T^{1/4+3\epsilon }$, then on the event $\mathcal E$ we have that $\bigcup _{s\in [t,t+T^{1/4+3\epsilon }]} \mathcal L_{s}'$ is contained in a ball of radius $CT^{1/4+3\epsilon }$ because any blue path to $\mathcal L_{s}'$ will have to pass through $\mathcal L_t'$ so it cannot get far within time $T^{1/4+3\epsilon }$. Note that $\mathcal E$ is an event in the original environment and not in the perturbed environment, but passage times do not change by much between the environments (by more than factor $2$ say). We obtain that in any time interval of length $T^{1/4+3\epsilon }$ we can decrease at most $C(T^{1/4+3\epsilon })^2$ many edge weights. We can cover $[0,T]$ with $\lceil  T^{3/4-3\epsilon } \rceil $ such intervals, which completes the proof of the claim.
\end{proof}

\begin{proof}[Proof of Lemma~\ref{lem:TV}] Since the total variation distance between two variables does not depend on whether they are coupled or not (and if so, in what way), we will use a different coupling of $P$ and $P'$. 
  
  To do so, we use two independent sequences of i.i.d.\ $\exp(1)$ variables, $Y_1,Y_2,\dots $ and $Z_1,Z_2,\dots $, thinking of the latter as a ``reservoir''. 
  To sample the original process $P$ the idea is to use the random variables in the reservoir for all the blue waiting times for which the blue boundary is thin. More precisely, suppose that at time $t$ a new waiting time along an edge needs to be sampled. Namely, at time $t$ a vertex $v$ in this edge becomes infected and its neighbor along the edge is healthy. In this case, if $t\le T$, $v\in \mathcal L_t$ and $L_t\le T^{1/4+3\epsilon }$ then we use the variables in the reservoir $Z_1,Z_2,\dots $ to sample the times it takes $v$ to infect its healthy neighbors. In any other case (i.e. $t>T$ or $v\notin \mathcal L_t$ or $L_t >T^{1/4+3\epsilon } $) we use the sequence $Y_1,Y_2,\dots $ to sample the infection times.

    To sample the perturbed process $P'$ we define the perturbed reservoir $Z_1',Z_2',\dots $ by $$Z_i'=\frac{Z_i}{1+\delta }.$$
    The process $P'$ is defined in the same way as $P$ with the same sequence $Y_1,Y_2,\dots $ and with the perturbed reservoir $Z_1',Z_2',\dots $. This is clearly a valid way of sampling both $P$ and $P'$. (The coupling we get from this process would be unusable for the previous part of the argument because, on the first time when one of them takes a variable from the reservoir and the other does not, they desynchronise --- each one takes values from $Y$ or $Z$ in order --- and diverge completely. Nevertheless, as we will see shortly, this coupling is useful for estimating the total variation distance).

    To bound the total variation distance between $P$ and $P'$, we define an auxiliary process $\Tilde{P}$ as follows. 
    Let $n=\lfloor T^{5/4+4\epsilon }\rfloor $ and let $\tilde Z_1,\tilde Z_2,\dots $ defined by 
\begin{equation}
    \tilde Z_i=\begin{cases}
     Z_i' \quad i\le n \\
     Z_i \quad i>n
    \end{cases}
\end{equation}
    The process $\tilde P$ is defined exactly as $P,P'$ but it uses $\tilde Z$ instead of $Z'$.
    
    By Lemma~\ref{lem:mermin} we have that $d_{\rm TV} (P,\tilde{P})\le d_{\rm TV}(Z,\tilde Z)\le \delta \sqrt{n} \le CT^{-\epsilon }$. Moreover, we have that  $\{N'\le n\} \subseteq \{P'=\tilde P\}$ and therefore by Claim~\ref{claim:k} we have $d_{\rm TV}(\tilde P, P')\le \mathbb P (N' \ge n)\le Ce^{-cT^c}$. This completes the lemma using the triangle inequality for the total variation distance.
\end{proof}

\section{Limit shape geometry: proof of Proposition~\ref{prop:limit}}\label{sec:limit}

The basic idea behind Proposition~\ref{prop:limit} is that, once a sufficiently large blue interval is created, the law of large numbers starts kicking in, it advances essentially deterministically and will soon overtake the red due to its faster speed. 
The main obstacle is that `large' here means $t^{1/4}$, which is too little for a straightforward application of the fluctuation bounds (\ref{non-random})-(\ref{random}). This will be handled by considering a second-order approximation $\tilde B_t$, the main object of this section. Recall that $B_t$ is the FPP ball coupled in the natural coupling with our $(W_t,S_t)$ process, and the comparisons (\ref{eq:WBWS}).

Throughout this section we assume that $T$ is sufficiently large and fix $t\in [T^{1/4},T]$. Our goal is to show that if the boundary at time $\approx t$ is larger than $T^{1/4+3\epsilon }$ then the strong type will encircle the weak one with high probability.

\begin{figure}
    \centering
        \begin{overpic}[width=0.6\linewidth]{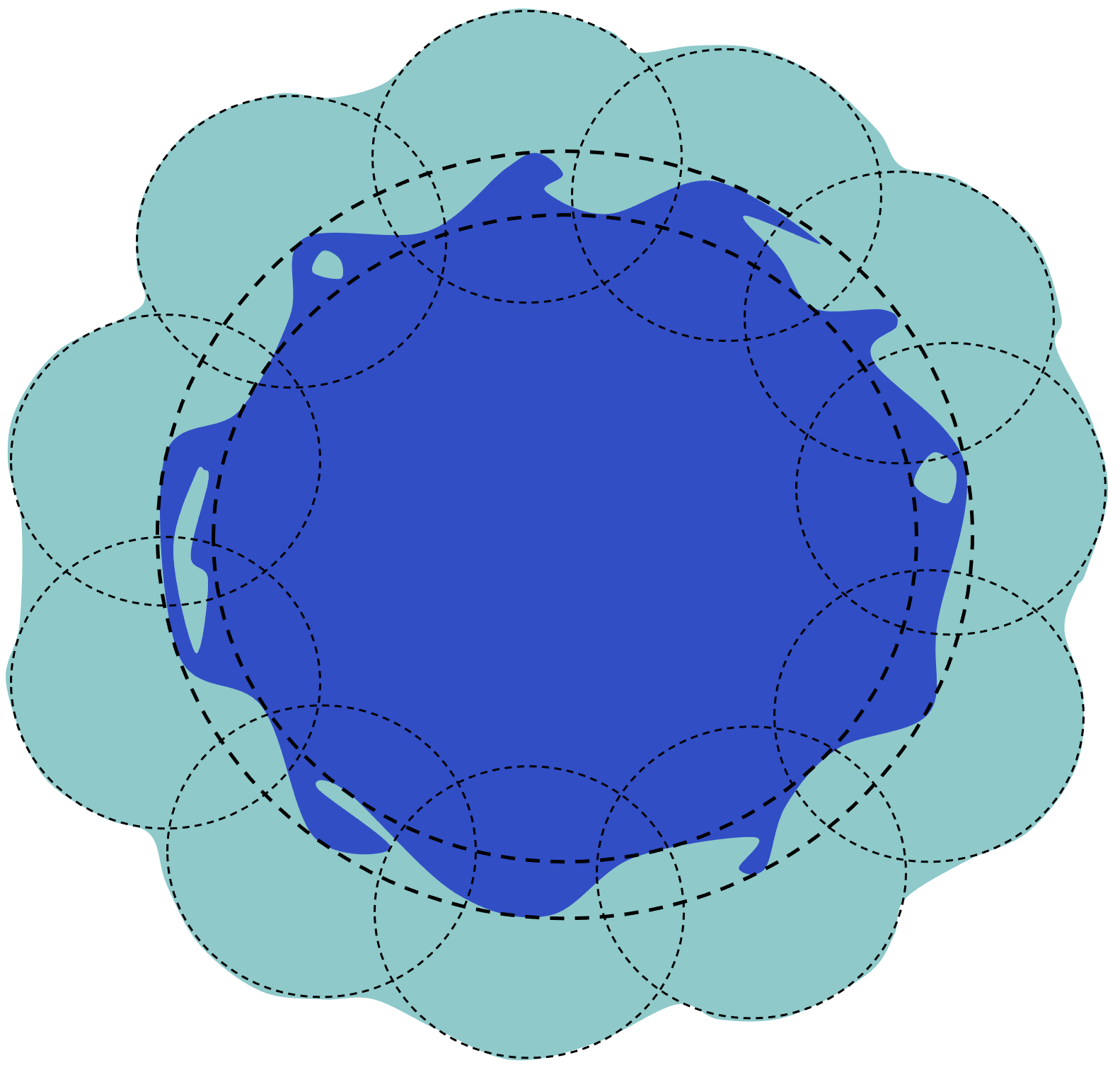}
        \put(47,57){\ensuremath{B_{t-t^{1/2+2\epsilon }}}}
        \put(83,69.5){\ensuremath{\tilde B_t}}
    \end{overpic}
    \caption{The second order approximation. While the balls $B_t$ and $B_{t-t^{1/2+2\epsilon }}$ may a priori contain holes and tunnels, the second order approximation $\tilde B_t$ is much more regular.}
    \label{fig:second}
\end{figure}

We will need the second order approximation of the ball $B_t$ (see Figure~\ref{fig:second}). Define
\begin{equation}\label{eq:second}
    \tilde B_t:= B_{t-t^{1/2+2\epsilon }} + t^{1/2+2\epsilon }\mathcal B. 
\end{equation}

$\tilde B_t$ is a better approximation to $B_t$ than simply $t\mathcal B$. The Hausdorff distance between $B_t$ and $\tilde B_t$ will be of order $t^{1/4+2\epsilon }$ (see \eqref{eq:3} and \eqref{eq:4} below for a more accurate statement).

Let 
\begin{equation}
    \delta :=\min \Big( \frac{\lambda -1}{8} \ , \ \frac{1}{20}\Big).
\end{equation}
($\lambda$ being, as usual, the rate of the faster, blue process $S_t$).

\begin{lem}\label{lem:limitshape}
 If the event  $\mathcal E$ holds, then for any $z\in \partial \tilde B_t$, the escaping ray $z+s\hat z$ with $\hat z=z/\mu (z)$ satisfies 
    \begin{equation}\label{eq:1}
      d_{\mu }(z+s\hat z, \tilde B_t ) \ge (1-\delta )s, \quad \text{ for all } s\ge 0.
    \end{equation}
\end{lem}

\begin{proof}
    Let $R=t^{1/2+2\epsilon }$ and $r=t^{1/2+\epsilon }$. On the event $\mathcal E$ we have that 
    \begin{equation}\label{eq:geom}
        (t-R-r) \mathcal B \subseteq B_{t-R} \subseteq (t-R+r) \mathcal B.
    \end{equation}
    This is the only use of $\mathcal E$ in this lemma, it is in effect a geometric lemma that holds for any set that satisfies (\ref{eq:geom}).
    
    It suffices to prove that for any $s>0$ and $y\in B_{t-R}$ we have that
    \begin{equation}\label{eq:2}
        \mu (z+s\hat z-y)\ge R+(1-\delta )s.
    \end{equation}
    Indeed, then for any $x\in y+R\mathcal B$ we have that $\mu (z+s\hat z-x)\ge \mu (z+s\hat z-y)-R\ge (1-\delta )s$, establishing \eqref{eq:1}.

    Fix $s>0$ and $y\in B_{t-R}$. Let $v=z+s\hat z-y$ and $\hat v=v/\mu (v)$. We may assume that $\mu (v)\le R+s$ as otherwise \eqref{eq:2} clearly holds. We start by proving that
    \begin{equation}\label{eq:eta}
        \mu \Big( \frac{\hat v+\hat z}{2} \Big) =\frac{\mu \big( v+\mu (v)\hat z \big)}{2\mu (v)} \ge 1-\frac{2r}{R} .
    \end{equation}
    Our assumption that $\mu(v)\le R+s$ bounds the denominator. To bound the numerator we use that $\mu (y)\le t-R+r$ and $\mu (z)\ge t-r$. These two give that $\mu (v)\ge \mu (z+s\hat z)-\mu (y)\ge R+s-2r$ and we obtain
    \begin{equation}
        \mu \big( v+\mu (v)\hat z \big) =\mu(z+\hat z(s+\mu(v))-y)\ge \mu (z)+s+\mu (v) -\mu (y) \ge 2R+2s-4r. 
    \end{equation}
    This completes the proof \eqref{eq:eta}. By Remark~\ref{rem:1} and \eqref{eq:eta}  we have that $\mu (\hat v-\hat z) \le \delta $ for large $t$ and therefore 
    \begin{equation}
        R\le \mu (z-y) =\mu (v-s\hat z) \le \mu (v-s\hat v) +s\mu(\hat v-\hat z) \le \mu (v)-s+\delta s, 
    \end{equation}
    which finishes the proof of \eqref{eq:2} and thus of Lemma \ref{lem:limitshape}.
\end{proof}

\begin{figure}
    \centering
        \begin{overpic}[width=0.6\linewidth]{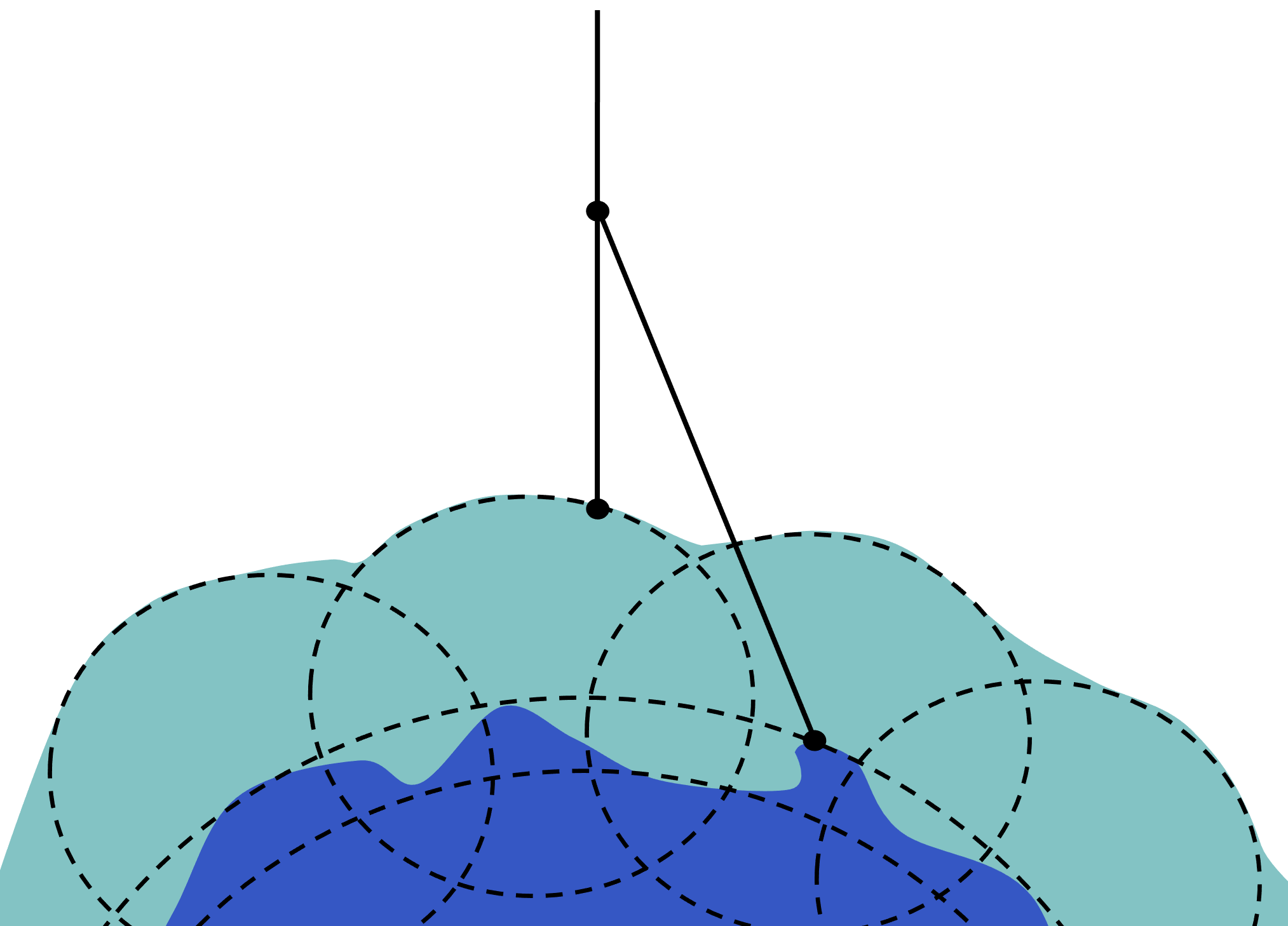}
        \put(48,55){\ensuremath{z+s\hat z}}
        \put(43.4,33.7){\ensuremath{z}}
        \put(64.3,15){\ensuremath{y}}
    \end{overpic}
    \caption{Proof of Lemma~\ref{lem:limitshape}.}
    \label{fig:7}
\end{figure}

Consider the times 
\begin{equation}
    t^-:= t-3t^{1/4+2\epsilon } \quad \text{and} \quad t^+:=t+3t^{1/4+2\epsilon }.
\end{equation}
With high probability when $s\le t^-$ we have $B_s\subseteq \tilde B_t$ and when $s\ge t^+$ we have $B_s\supseteq \tilde B_t$. In fact, we claim that on  $\mathcal E$ we have
\begin{equation}\label{eq:3}
      \tilde B_t \subseteq B_{t-t^{1/2+2\epsilon }} + (t^{1/2+2\epsilon } +2t^{1/4+2\epsilon})\mathcal B \subseteq  B_{t^+} \stackrel{\smash{(\ref{eq:WBWS})}}{\subseteq} W_{t^+}\cup S_{t^+}
\end{equation}
and
\begin{equation}\label{eq:4}
    W_{t^-} \stackrel{\smash{(\ref{eq:WBWS})}}{\subseteq} B_{t^-} \subseteq B_{t-t^{1/2+2\epsilon }} + (t^{1/2+2\epsilon }-2t^{1/4+2\epsilon})\mathcal B \subseteq \tilde B_t.
\end{equation}
The second inclusion in \eqref{eq:3} follows since on $\mathcal E$ every vertex at $\mu $ distance at most $t^{1/2+2\epsilon }+2t^{1/4+2\epsilon}$ from   a vertex in $B_{t-t^{1/2+2\epsilon }}$ will have passage time at most $t^{1/2+2\epsilon }+3t^{1/4+2\epsilon}$ to that vertex.

To prove the second inclusion in \eqref{eq:4} take $y\in B_{t^-}$ and consider the geodesic $\gamma $ connecting $0$ and $y$ in the FPP environment. Let $x$ be the first point of $\gamma $ outside of $B_{t-t^{1/2+2\epsilon }}$. We have that $T(x,y)\le t^--(t-t^{1/2+2\epsilon })= t^{1/2+2\epsilon }-3t^{1/4+2\epsilon }$ and therefore using that $x$ is a neighbor of a vertex in $B_{t-t^{1/2+2\epsilon }}$ we get that on $\mathcal E$
\begin{equation}
    d_{\mu }(y,B_{t-t^{1/2+2\epsilon }}) \le \mu (x-y) +d_{\mu }(x,B_{t-t^{1/2+2\epsilon }}) \le t^{1/2+2\epsilon } -2t^{1/4+2\epsilon }. 
\end{equation}

\begin{figure}
    \centering
        \begin{overpic}[width=0.7\linewidth]{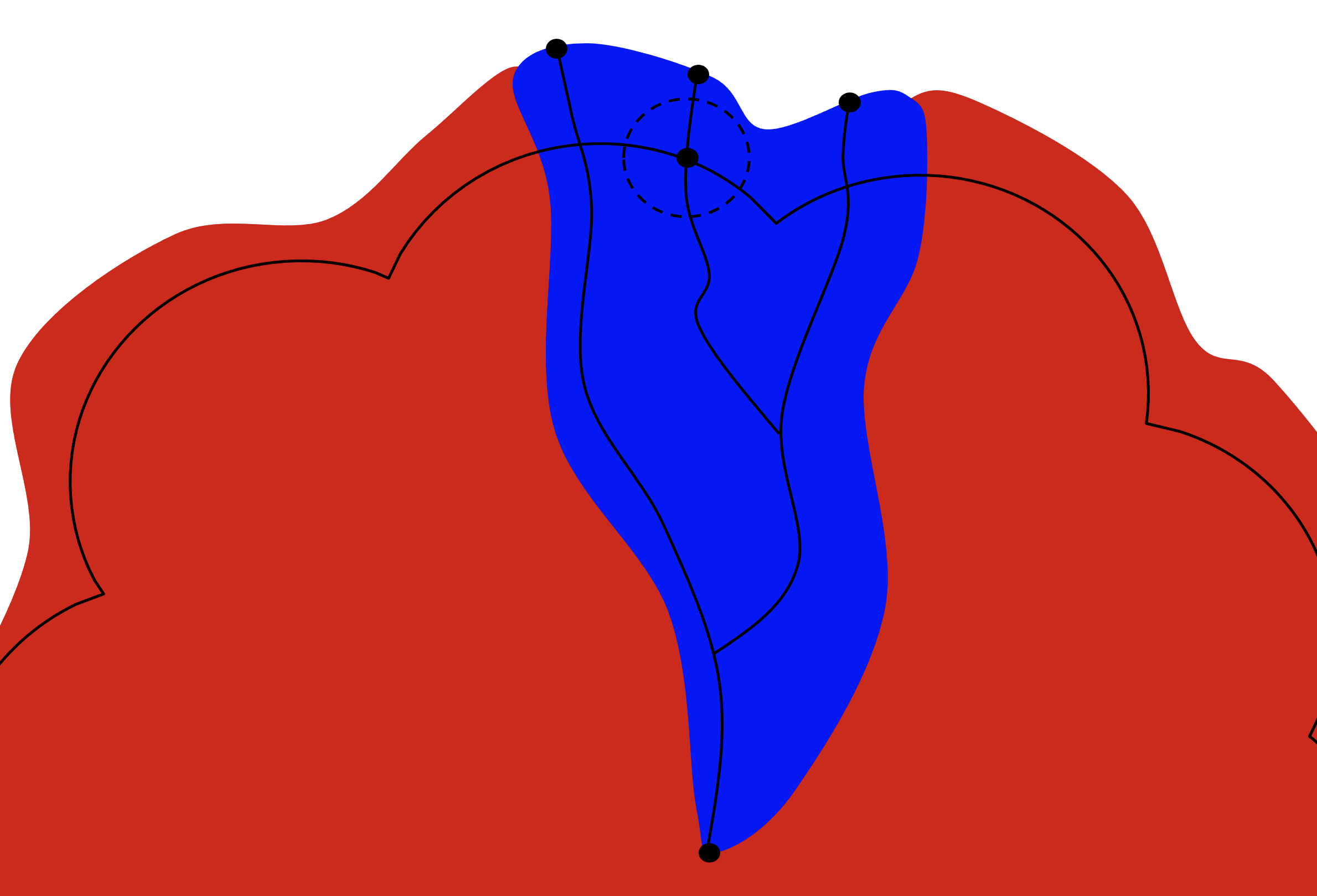}
        \put(50,1.8){\ensuremath{v_0}}
        \put(40,65.5){\ensuremath{u}}
        \put(52,63.7){\ensuremath{w}}
        \put(62,61.3){\ensuremath{v}}
        \put(49.5,54){\ensuremath{z}}
        \put(44.3,41){\ensuremath{\gamma _u}}
        \put(53.2,43){\ensuremath{\gamma _w}}
        \put(61.9,42.5){\ensuremath{\gamma _v}}
        \put(19,44.5){\ensuremath{\partial \tilde B_t}}
    \end{overpic}
    \caption{Proof of Lemma~\ref{lem:ball}. }
    \label{fig:8}
\end{figure}

\begin{lem}[big blue ball]\label{lem:ball}
  Suppose that $L_{t^+}\ge t^{1/4+3\epsilon }$ 
  and that $\mathcal E$ holds. Then there exists $z\in \partial \tilde B_t$ such that $z+t^{1/4+2\epsilon }\mathcal B\subseteq S_\infty $.
\end{lem}

\begin{proof}
 Suppose that $L_{t^+}\ge t^{1/4+3\epsilon }$ and that $\mathcal E$ holds and define
\[
    \tau := \min \big\{  s\ge 0  : S_s \setminus \tilde B_t \neq \emptyset  \big\}.
\]
i.e.~the first time blue exits $\tilde B_t$. Note that blue grows linearly faster and so $\tau$ can be much smaller than $t$.

First, suppose that $\tau \le t-5t^{1/4+2\epsilon }$ (i.e.~blue is significantly ahead of $B_t$) and let $z$ be the first point of $\partial \tilde B_t$ infected by blue (it exists since $\mathcal L_{t^+}\neq \emptyset $ and $\mathcal L_{t^+}\cap \tilde B_t=\emptyset$ by \eqref{eq:3}). In this case, red will be very far at time $\tau $ and blue has plenty of time to capture a ball around $z$ of radius $t^{1/4+2\epsilon }$. Indeed, once $z$ is infected at time $\tau $, by time $\tau ^+:=\tau +2t^{1/4+2\epsilon }\le t^-$ the entire ball $z+t^{1/4+2\epsilon }\mathcal B$ will be infected (by either red or blue). However, by \eqref{eq:4} we have that $d_{\mu }(z,W_{\tau ^+}) \ge d_{\mu }(z,W_{t^-}) \ge 2t^{1/4+2\epsilon }$ which implies that $z+t^{1/4+2\epsilon }\mathcal B \subseteq S_{\tau ^+}\subseteq S_{\infty }$.

Next, suppose that $\tau \ge t-5t^{1/4+2\epsilon }$. Since $L_{t^+}\ge t^{1/4+3\epsilon }$, there are $u,v\in \mathcal L_{t^+}$ such that $\|u-v\|\ge t^{1/4+3\epsilon }$. Let $\gamma _u,\gamma _v$ be the blue paths of infections from $v_0$ to $u,v$ respectively (see Figure~\ref{fig:8}). By planarity, the union of the paths $\gamma _u \cup \gamma _v$ separates the infected region $S_{t^+}\cup W_{t^+}$ into two disjoint parts, and we let $S_{u,v}$ be the part that does not contain the origin together with $\gamma _u\cup \gamma _v$ (more precisely, $S_{u,v}$ are all vertices in $S_{t^+}\cup W_{t^+}$ such that any path in $S_{t^+}\cup W_{t^+}$ connecting them to the origin intersects $\gamma _u\cup \gamma _v$). The red infection starting from the origin cannot enter this region and so it must be blue $S_{u,v}\subseteq S_{t^+}$. We would like to find a ball inside $S_{u,v}$. Let $w\in \mathcal L_{t^+}\cap S_{u,v}$ on the part of the blue boundary between $u$ and $v$ which is at distance at least $\frac{1}{3}t^{1/4+3\epsilon }$ from both (see Figure~\ref{fig:8}). Note that $u,v,w\notin \tilde B_t$ since $\mathcal L_{t^+} \cap \tilde B_t=\emptyset$ by \eqref{eq:3}. Let $\gamma _w$ be the blue path of infections from $v_0$ to $w$ and let $z$ be the last intersection of $\gamma _w$ with $\partial \tilde B_t$. Since $\tau \ge t-5t^{1/4+2\epsilon }$ it follows that on $\mathcal E$ we have $\|z-w\|\le Ct^{1/4+2\epsilon }$. By \eqref{eq:3} we have that $z+t^{1/4+2\epsilon }\mathcal B\subseteq S_{t^+}\cup W_{t^+}$ and we claim that in fact on $\mathcal E$ we have $z+t^{1/4+2\epsilon }\mathcal B\subseteq S_{t^+}$. Indeed, $\gamma _w$ will be sandwiched between $\gamma _u$ and $\gamma _v$ and will be contained in $S_{u,v}\subseteq S_{t^+}$ (see Figure~\ref{fig:8}). Hence, the only way in which $z+t^{1/4+2\epsilon }\mathcal B$ can contain a red point, is if either $\gamma _u$ or $\gamma _v$ intersect $z+t^{1/4+2\epsilon }\mathcal B$. However $\gamma_u,\gamma _v$ cannot enter this ball before time $\tau -2t^{1/4+2\epsilon }$ since otherwise on $\mathcal E$,  $z$ will be infected before time $\tau $ contradicting the definition of $\tau $. Thus, using our assumption that $\tau \ge t-5t^{1/4+2\epsilon }$, this means that $\gamma _u,\gamma _v$ had only $Ct^{1/4+2\epsilon }$ time to reach from $z+t^{1/4+2\epsilon }\mathcal B$ to their respective endpoints $u,v$, which is impossible on $\mathcal E$, since they had to travel distance $ct^{1/4+3\epsilon }$.
\end{proof}

\begin{lem}[Escaping]\label{lem:escape}
    Suppose that $L_{t^+}\ge t^{1/4+3\epsilon }$ and that $\mathcal E$ holds. Then, there exists $u\in S_{2t}$ with $\mu (u) \ge (2+3\delta )t$.
\end{lem}

\begin{proof}
Let $z\in \partial \tilde B_t$ be the point from Lemma~\ref{lem:ball} with $z+t^{1/4+2\epsilon } \mathcal B \subseteq S_{\infty }$. Let $\tau _0$ be the infection time of $z$. On $\mathcal E$, by the same arguments as in \eqref{eq:3} we have that $\tau _0\le t+t^{1/4+3\epsilon /2}$ and by the same arguments as in \eqref{eq:4} we have that
\begin{equation}\label{eq:10}
    W_{\tau _0} \subseteq B_{\tau _0} \subseteq  \tilde B_t + 2t^{1/4+3\epsilon /2}\mathcal B.
\end{equation}

Starting from $z_0=z$ we define a sequence of points by 
$$z_i:=z+it^{1/4}\hat z, \quad i\ge 0$$ 
(as before, $\hat z=z/\mu(z)$) and let $\tau _i$ be the first time $z_i$ is infected. We prove using induction on $i\le t$ that 
\begin{equation}\label{eq:induction}
    \tau _i \le \tau _0+ (1-5\delta )it^{1/4} \quad \text{ and that }\quad z_i\in S_{\tau _i}.
\end{equation}
This clearly holds when $i=0$. Next, we assume that \eqref{eq:induction} holds and prove it at step $i+1$.

Let $T_i=T(z_i,z_{i+1})$ be the passage time from $z_i$ to $z_{i+1}$ in the environment $\{X_e\}_{e\in E(\mathbb Z^2)}$. Let $\gamma _i$ be the corresponding geodesic from $z_i$ to $z_{i+1}$. On $\mathcal E$ we have that $T_i\le (1+\delta )t^{1/4}$ and that $\gamma _i \subseteq z_i+ 2t^{1/4}\mathcal B$ (as otherwise $\gamma _i$ contains a point at distance $2t^{1/4}$ from $z_i$ whose passage time to $z_i$ is too small). Thus, for the induction in \eqref{eq:induction} it suffices to prove that  
\begin{equation}\label{eq:9}
    d_\mu (z_i,W_{\tau _{i+1}}) \ge 3t^{1/4}. 
\end{equation}
Indeed, if \eqref{eq:9} holds  then all of $\gamma _i$ is blue and therefore $\tau _{i+1}-\tau _i \le T_i/\lambda \le (1-5\delta )t^{1/4}$ completing the induction.

We turn to prove \eqref{eq:9}. When $i\le \frac{1}{2}t^{2\epsilon }$ this estimate follows since $z+t^{1/4+2\epsilon } \mathcal B \subseteq S_{\infty }$. Suppose next that $i\ge \frac{1}{2}t^{2\epsilon }$. Note that $\tau _{i+1}-\tau _i\le T_i\le (1+\delta )t^{1/4}$ and therefore by the induction hypothesis $\tau _{i+1}-\tau _0\le (1-4\delta )it^{1/4}$. This means that 
\begin{equation}
    W_{\tau _{i+1}} \subseteq W_{\tau _0} + (1-3\delta )it^{1/4} \mathcal B \subseteq \tilde B_t + (1-2\delta )it^{1/4} \mathcal B,
\end{equation}
where in the last inclusion we used \eqref{eq:10}. This finishes the proof of \eqref{eq:9} since by Lemma~\ref{lem:limitshape} we have that $d_{\mu } (z_i,\tilde B_t) \ge (1-\delta )it^{1/4}$.

Once the induction in \eqref{eq:induction} is completed, the lemma follows by taking $i=\lfloor t^{3/4}/(1-5\delta )- t^{1/2+3\epsilon /2}\rfloor $. For such $i$ we have that $\tau _i \le \tau _0 +(1-5\delta )it^{1/4} \le 2t$ and moreover $\mu (z_i) = \mu (z)+ it^{1/4} \ge (2+3\delta )t$, as needed.
\end{proof}

\begin{figure}
    \centering
        \begin{overpic}[width=0.6\linewidth]{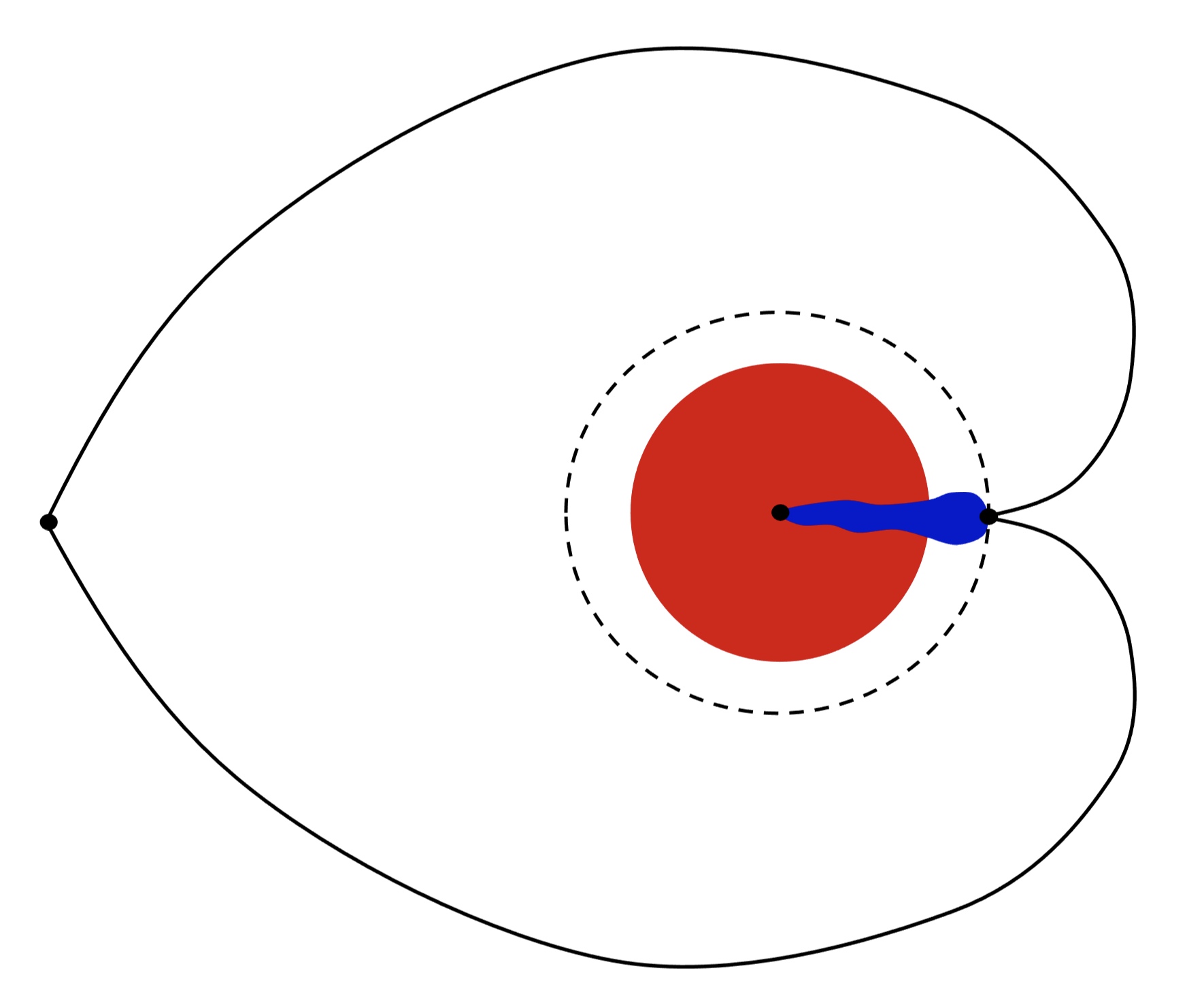}
        \put(82.7,37.3){\ensuremath{u}}
        \put(88,70){\ensuremath{\gamma _1}}
        \put(88,11.2){\ensuremath{\gamma _2}}
        \put(63,37.5){\ensuremath{0}}
    \end{overpic}
    \caption{Proof of Lemma~\ref{lem:kill}. }
    \label{fig:6}
\end{figure}

\begin{lem}[Killing red]\label{lem:kill}
   Suppose that $L_{t^+}\ge t^{1/4+3\epsilon }$ and that $\mathcal E$ holds. Then  $|W_{\infty }|<\infty $.
\end{lem}

\begin{proof}
  By Lemma~\ref{lem:escape}, there exists $u\in S_{2t}$ with $\mu (u)\ge (2+3\delta )t$. Moreover, on the event $\mathcal E$ we have that $W_{2t}\subseteq B_{2t} \subseteq (2+\delta )t\mathcal B$. We will show that when the strong type is that far from the weak one it will quickly encircle it.
    
    Let $R=\mu (u)$. We claim that there are two curves $\gamma _1:[0,\ell _1]\to \mathbb R ^2$ and $\gamma _2:[0,\ell _2]\to \mathbb R ^2$ such that the following holds:
  \begin{enumerate}
      \item
      The curves start from $\gamma _1(0)=\gamma _2(0)=u$.
      \item 
      The curves are in natural parametrization with respect to $\mu$. Namely, $\gamma _i'(s)\in \partial \mathcal B$ for $i=1,2$ and $s\in [0,\ell _i]$. This implies that $\mu (\gamma _i(s)-\gamma _i(s'))\le |s-s'|$.
      \item 
      For $i=1,2$ and any $s\in [0,\ell _i]$ we have that $\mu (\gamma _i (s))\ge R+(1-\delta )s$.
      \item 
      The curves have linear length: $\ell _i\le Ct$.
      \item 
      They have the same endpoint $\gamma _1 (\ell _1)=\gamma _2(\ell _2)$.
      \item 
      Together the curves circle the origin. Namely, $\gamma _1\cup \gamma _2$ separates $0$ from infinity.
  \end{enumerate}  
See Figure~\ref{fig:6}. To prove the existence of these curves, let $\theta _1,\theta _2$ be two curves in natural parametrization starting at $\hat u=u/\mu (u)$ and going clockwise and counter clockwise around $\partial \mathcal B$ until reaching $-\hat u$. Define the curves 
\begin{equation}
    \Gamma _i(s):= R e^s \theta _i(\delta s).
\end{equation}
We have that $\mu (\Gamma _i(s))=Re^s$. Moreover, $\Gamma _i'(s)=Re^s(\theta _i(\delta s)+\delta \theta _i'(\delta s))$ and therefore $\mu (\Gamma _i'(s)) \le Re^s(1+\delta )$. This gives that $\mu (\Gamma _i(s))-R \ge \frac{1}{1+\delta }\int _0^s \mu (\Gamma _i'(r))dr$, but this integral is precisely the length of the curve $\Gamma _i$ up to time $s$. Thus, letting $\gamma _i$ be the natural parametrization of $\Gamma _i$ we have that $\mu (\gamma _i(u))-R\ge \frac{u}{1+\delta } \ge  (1-\delta )u$. The other requirements on $\gamma _1,\gamma _2$ clearly hold.

On the event $\mathcal E$, the strong type will infect both curves before the weak one will have a chance to get close and interrupt it. This is made rigorous by taking a sequence of points along $\gamma _i$ as in the proof of Lemma~\ref{lem:escape}. The details are omitted.
\end{proof}

We can now prove Proposition~\ref{prop:limit}.

\begin{proof}[Proof of Proposition~\ref{prop:limit}]
    By Lemma~\ref{lem:E}, it suffices to prove that for large $T$ we have
\begin{equation}\label{eq:inclusion}
    \mathcal E \subseteq \big\{ |W_\infty |<\infty  \big\} \cup \big\{ \forall t\le T, \ L_t < T^{1/4+3\epsilon} \big\}.
\end{equation}
First, we note that if $t\le 2T^{1/4}$, then on $\mathcal E$ we have that $S_t\subseteq 3\lambda T^{1/4}\mathcal B$ which implies that $L_t< T^{1/4+3\epsilon }$.

Next, if there exists a time $s\in [2T^{1/4},T]$ such that $L_s\ge T^{1/4+3\epsilon }$ then there exists a time $t\in[T^{1/4},T]$ such that $s=t^+$ and then by Lemma~\ref{lem:kill} we have that $|W_\infty |<\infty $. This finishes the proof of \eqref{eq:inclusion}.
\end{proof}

\bibliographystyle{abbrv}
\bibliography{bibliography}

\end{document}